\documentclass[a4paper,10pt]{amsart}
\usepackage[english]{certus}
\usepackage{a4wide}
\usepackage[textsize=footnotesize,linecolor=lime,backgroundcolor=lime]{todonotes}
\usepackage{enumitem}

\makeatletter
\newcommand{\mylabel}[2]{#2\def\@currentlabel{#2}\label{#1}}
\makeatother

\theoremstyle{definition}
\title{Maximal amenable subgroups of arithmetic groups}
\author{Vadim Alekseev}
\address{Vadim Alekseev, Technische Universit\"{a}t Dresden, Fakult\"{a}t Mathematik, Institut f\"{u}r Geometrie, 01062 Dresden, Germany}
\email{vadim.alekseev@tu-dresden.de}

\author{Alessandro Carderi}
\address{Alessandro Carderi, Institut für Algebra und Geometrie, Karlsruhe Institute of Technology, 76128 Karlsruhe, Germany}
\email{alessandro.carderi@kit.edu}


\begin{document}
\onehalfspace

\begin{abstract}
By classifying $S$-maximal amenable subgroups of algebraic groups over a global field of characteristic zero, we obtain a complete classification of maximal amenable subgroups up to commensurability in the respective arithmetic groups. Futhermore, we prove that these commensurably maximal amenable subgroups are singular and therefore give rise to maximal amenable von Neumann subalgebras.
\end{abstract}

\maketitle

\section{Introduction}

In this article we are studying amenable subgroups of arithmetic groups and their ambient locally compact groups. By the Tits alternative \cite{TitsFree1972}, all amenable subgroups of arithmetic groups are virtually solvable and hence structurally easy to understand; similarly, amenable real Lie groups are composed of solvable and compact real Lie groups and hence again structurally easy to handle. However, understanding maximal amenable subgroups of a given group is a much harder task. For real Lie groups this has been accomplished by Moore \cite{MooreAmenable1979}. The primary goal of this article is to understand maximal amenable subgroups of arithmetic groups; since the latter ultimately come from algebraic groups, our analysis heavily relies on structure theory of algebraic groups wherefore our results can also be understood as an extension of Moore's work. 

Classifying all maximal amenable subgroups of arithmetic groups is a challenging task due to presence of finite maximal amenable subgroups, for instance, $\mb Z/6$ in $\mf{SL}_2(\mb Z)\cong \mb Z/6\ast_{\mb Z/2} \mb Z/4$. However, the situation becomes much more tractable if one works up to finite index sub- and overgroups, thus considering maximal amenable commensurability classes instead of maximal amenable subgroups; every such class contains a maximal amenable subgroup. For example, in the realm of hyperbolic groups the work of Gromov \cite{GromovHyperbolic1987} yields complete classification of infinite maximal amenable subgroups: they are stabilisers of couples of points on the Gromov boundary. In particular, their commensurability classes are maximal amenable.

One of our main results is the complete classification of commensurably maximal amenable subgroups in arithmetic groups over a field of characteristic zero. They correspond to algebraic subgroups satisfying a maximal amenability condition with respect to a given set $S$ of places of the ground field. Surprisingly, these algebraic subgroups can be bigger than the Zariski closure of the maximal amenable subgroup, so that establishing the correspondence requires some delicate considerations. The natural counterpart to commensurability in the setting of algebraic groups considered over a given set of places is working up to $S$-cocompact sub- and overgroups, meaning that over every place $v\in S$ the group becomes cocompact. With this idea, the relevant condition on the algebraic subgroup $\mf H\leqslant \mf G$ turns out to be ``$S$-cocompact maximal amenability'', and our classification result describes $S$-cocompactly maximal amenable subgroups $\mf H$ of a given algebraic group $\mf G$. It turns out that upon fixing the unipotent radical $\mf U$ of such a subgroup the quotient has an $S$-cocompact subtorus. This torus has no $S$-cocompact subgroups and is $S$-cocompact in its centralizer. Together with a rank condition which ensures that the torus normalizes $\mf U$ and not a bigger unipotent subgroup, we get the definition of an $S$-ample torus (see Definition \ref{def:s-ample}) which allows us to formulate the classification result: 

\begin{ThmA}[Theorem \ref{thm:maximal-amenable-algebraic-subgroups}]
Let $K$ be a global field of characteristic zero, $\mf G$ be a connected $K$-group and  $S\subseteq \ms R$ be a set of valuations. 
Every $S$-cocompactly maximal amenable subgroup $\mf H\leqslant \mf G$ is $S$-cocompactly equivalent to a subgroup $\mf H_{\mf T}\leqslant \mf H$ of the form
\[
\mf H_{\mf T} \coloneqq \mf T \ltimes \mf U,
\]
where $\mf U$ is the unipotent radical of the $K$-parabolic subgroup $\mf P\coloneqq \mc N_{\mf G}(\mf U)$ of $\mf G$ with a Levi decomposition $\mf P = \mf L\ltimes \mf U$, and $\mf T\leqslant \mf L$ is an $S$-ample torus. Conversely, the normalizer of every subgroup $\mf H_{\mf T}$ as above is $S$-cocompactly maximal amenable. 
\end{ThmA}

By taking the $K(S)$-points, we obtain the corresponding classification of commensurably maximal subgroups in arithmetic groups. Let us remark that we don't assume $S$ to be finite; in particular, our result also cover commensurably maximal amenable subgroups of $\mf G(K)$.

\begin{ThmA}[see Proposition \ref{prop: every CMA is H(K(S))}]
  A subgroup $\Lambda<\mf G(K(S))$ is commensurably maximal amenable if and only if there is an $S$-cocompactly maximal amenable subgroup $\mf H<\mf G$ such that $\Lambda=\mf H(K(S))$.
\end{ThmA}

Furthermore, we prove that all commensurably maximal amenable subgroups verify a stronger condition: they are singular in the ambient arithmetic groups.

\begin{Def*}
    An amenable subgroup $\Lambda<\Gamma$ is singular if there is a continuous action $\Gamma\curvearrowright X$ on a compact space such that for any $\Lambda$-invariant probability measure $\mu$ and for all $\gamma\in \Gamma\setminus \Lambda$, the probability measures $\mu$ and $\gamma_*\mu$ are singular.
\end{Def*}

The notion of singularity comes from the context of von Neumann algebras, it was introduced in \cite{BoutonnetMaximal2015}, where it was shown that all singular subgroups give rise to maximal amenable subalgebras. Currently, this is the main source of maximal amenable subalgebras.

In hyperbolic groups every commensurably maximal amenable subgroup is singular with respect to the action on the Gromov boundary \cite[Proposition 2.1]{BoutonnetMaximal2015}. Other examples of singular subgroups were obtained for acylindrically hyperbolic groups in \cite[Theorem 3.2]{JiangSingular2020}. Beyond hyperbolicity, singularity was known for HNN extensions and amalgamated free products \cite[Proposition 2.4, Proposition 2.5]{BoutonnetMaximal2015} and for certain particular examples of maximal amenable subgroups of lattices in real Lie groups \cite[Proposition 2.6, Proposition 2.7]{BoutonnetMaximal2015}. Our results extend this to all commensurably maximal amenable subgroups in arithmetic groups over a field of characteristic zero, completely classifying singular subgroups in these groups.

\begin{ThmA}[Theorem \ref{ccma singular K}]
  Let $K$ be a global field of characteristic $0$ and let $S$ be a set of places of $K$ satisfying $\ms R_\infty\setminus \ms T(\mf G)\subseteq S$. Let $\mf G$ be a connected $K$-group and $\mf H$ be an $S$-cocompactly maximal amenable subgroup. Then $\Lambda \coloneqq \mf H(K(S))$ is a singular subgroup of $\Gamma \coloneqq \mf G(K(S))$. In particular, its von Neumann algebra $L\Lambda$ is maximal amenable in $L\Gamma$.
\end{ThmA} 

The space $X$ which is responsible for singularity of $\Lambda$ can be easily definied: it is just the product of the Furstenberg boundaries over places in $S$. However, the proof that we present here is much more complicated than in the previous cases: in our context the measure preserved by the amenable subgroup can be diffuse and the support of the translated measure is not necessarily disjoint, as was the case in all previous examples. This will force us to do a much finer analysis of the action on the Furstenberg boundary, heavily using the structure obtained in our classification result.

\vspace{3pt}

\paragraph{\textbf{Acknowledgments.}} The authors would like to thank Rémi Boutonnet, Yves de Cornulier, Cyril Houdayer and Andreas Thom for their interest in this work and helpful comments.  V.~A.\ acknowledges funding by the Deutsche Forschungsgemeinschaft (SPP 2026). A.~C.\ acknowledges funding by the Deutsche Forschungsgemeinschaft (DFG, German Research Foundation) – 281869850 (RTG 2229).

\section{Preface and notations}

In this text, we will follow the notations and conventions of \cite{MargulisDiscrete1991}. Throughout the paper, $K$ will denote a global field of characteristic $0$, that is, a finite extension of $\mb Q$. We denote by $\ms R$ the set of (equivalence classes of) the nonequivalent valuations of $K$ and by $\ms R_\infty \subset \ms R$ the set of Archimedean valuations. Given a valuation $v\in\ms R$, we will denote by $K_v$ the completion of $K$ with respect to $v$. The ring of adeles of $K$ will be denoted by $\mb A_K$, and given a subset $S\subseteq \ms R$, we denote by $K(S)$ the ring of $S$-integral elements of the field $K$. 
 
By a $K$-group, we mean an affine algebraic group defined over $K$. All our algebraic groups will be $K$-defined unless stated otherwise. 
Let $\mf H\leqslant \mf G$ be $K$-groups. We will use the following notations:
\begin{itemize}
  \item $\mf G^\circ$ will denote the connected component of the identity;
  \item $\mc R(\mf G)$ will denote the solvable radical of $\mf G$;
  \item $\mc R_u(\mf G)$ will denote the unipotent radical of $\mf G$;
  \item $\mc D(\mf G)$ will denote the derived subgroup of $\mf G$;
  \item $\mc Z(\mf G)$ will denote the center of $\mf G$;
  \item $\mc C_{\mf G}(\mf H)$ will denote the centralizer of $\mf H$ in $\mf G$;
  \item $\mc N_{\mf G}(\mf H)$ will denote the normalizer of $\mf H$ in $\mf G$;
  \item if $\mf T$ is a torus, $\mf T_a$ will denote its anisotropic part. 
\end{itemize}

Observe that all the above defined objects are $K$-groups \cite{BorelGroupes1965}, \cite{MargulisDiscrete1991}. 

A group $\mf G$ is \textit{reductive} if $\mc R_u(\mf G)$ is trivial and \textit{semisimple} if $\mc R(\mf G)$ is trivial. Note that $\mf G/\mc R_u(\mf G)$ is always reductive and that $\mf G/\mc R(G)$ is always semisimple. 
 Any $K$-group $\mf G$ has a reductive subgroup $\mf G^R$ called a Levi subgroup such that $\mf G\cong \mf G^R\ltimes \mc R_u(\mf G)$. This isomophism is called the \textit{Levi decomposition} of $\mf G$. Every reductive $K$-subgroup of $\mf G$ is conjugate to a subgroup of $\mf G^R$ by an element of $\mc R_u(\mf G)$ \cite[0.8]{BorelGroupes1965}.

Following \cite[I.3.2]{MargulisDiscrete1991}, given a subset $S\subseteq \ms R$, we denote by $\mf G_S$ the subgroup of $\mf G(\mb A_K)$ consisting of adeles whose $v$-components are $1$ for $v\not\in S$; it can also be interpreted as a restricted product of $\mf G(K_v)$ for $v\in S$. The group $\mf G_S$ is a locally compact topological group. If $S$ is finite, we get $\mf G_S= \prod_{v\in S}\mf G(K_v)$.

\section{$S$-maximal amenable algebraic subgroups}

In this section, we will fix a non-empty subset $S\subseteq \ms R$.

Let $\mf G$ be an algebraic group defined over $K$ and let $\mf H\leqslant\mf G$ be a $K$-defined subgroup. 

\begin{Def}
 We say that $\mf H$ is $S$-\textbf{amenable} if $\mf H(K_v)$ is amenable as locally compact group for every $v\in S$. 
\end{Def}

Remark that since an infinite direct sum of groups is dense in both the restricted and the unrestricted product, $\mf H$ is $S$-amenable if and only if $\prod_{v\in S} \mf H(K_v)$ is if and only if $\mf H_S$ is. 

\begin{Def}
   We say that $\mf H$ is $S$-\textbf{compact} if for every $v\in S$ the group $\mf H(K_v)$ is compact (with respect to the locally compact topology coming from $K_v$).
\end{Def}

Clearly $S$-compact subgroups are $S$-amenable.

\begin{Def}
An $S$-amenable subgroup $\mf H\leqslant\mf G$ is called $S$-\textbf{maximal amenable} if whenever $\mf H'\geqslant \mf H$ is a $K$-defined $S$-amenable group of $\mf G$, we have that $\mf H'=\mf H$. 
\end{Def}

First, we prove that every $S$-amenable subgroup $\mf H \leqslant \mf G$ is contained in an $S$-maximal amenable subgroup.

\begin{Prop}\label{prop: every SA contained in SMA}
  Let $\mf G$ be an algebraic group. For every $S$-amenable group $\mf H\leqslant \mf G$, there is an $S$-maximal amenable $\mf H'\leqslant \mf G$ containing $\mf H$.
\end{Prop}

In order to prove the above proposition, we will need the following consequence of the Jordan-Schur theorem. 

\begin{Lemma}\label{lemma:jordan schur}
Let $\mf G$ be an algebraic group and let $\mf N\normal \mf G$ be a normal subgroup. Suppose we have an increasing sequence of subgroups $\{\mf H_n\leqslant \mf G\}_n$ satisfying $\mf H_n^\circ=\mf N$. Then the Zariski closure of $\bigcup_n \mf H_n$ in $\mf G$ is an extension of a virtually abelian group $\mf A$ by $\mf N$.
\end{Lemma}
\begin{proof}
 Let $q\colon \mf G\to \mf G/\mf N$ be the natural quotient map. By assumption, the image $\Lambda_n \coloneqq q(\mf H_n)$  in $\mf G/\mf N$ is finite, hence the union $\Lambda \coloneqq \bigcup_n\Lambda_n$ is locally finite. Now, the Jordan--Schur theorem implies \cite[Corollary 12.1.12]{PassmanAlgebraic1977} that $\Lambda$ is virtually abelian, and therefore so is its Zariski closure in $\mf G/\mf N$. The claim follows.
\end{proof}

\begin{proof}[Proof of Proposition \ref{prop: every SA contained in SMA}]
  Let $\mf H$ be an $S$-amenable subgroup of $\mf G$. Clearly there is an $S$-amenable subgroup $\mf H_1\leqslant \mf G$ containing $\mf H$ of maximal dimension. Therefore $\mf H_1$ has finite index in any $S$-amenable subgroup of $\mf G$ containing it. Let $\{\mf H_n\}_n$ be an increasing chain of amenable subgroups containing $\mf H_1$. By assumption, the subgroup $\mf H_n^\circ=\mf H_1^\circ\eqqcolon \mf N$ is normal in $\mf G'\coloneqq \mc N_{\mf G}(\mf N)$. The group $\mf G'$ contains $\mf H_n$ for every $n$ and hence we can apply Lemma \ref{lemma:jordan schur} to obtain that the Zariski closure of $\bigcup_n \mf H_n$ is $S$-amenable. Since the dimension of $\mf H_1$ is already maximal, the sequence $\{\mf H_n\}_n$ stabilizes. 
\end{proof}

The aim of this section is to study $S$-maximal amenable subgroups of a given group $\mf G$. To do so, we perform several reduction steps. First of all, in view of the following lemma we may assume that $\mf G$ is reductive.

\begin{Lemma}\label{lem:ma-contains-unipotent-radical}
  Let $\mf G$ be a $K$-group. Then every $S$-maximal amenable subgroup $\mf H$ contains the unipotent radical of $\mf G$, that is $\mf H\geqslant \mc R_u(\mf G)$.
\end{Lemma}
\begin{proof}
  Assume that $\mf H\leqslant \mf G$ is $S$-amenable. Since $\mc R_u(\mf G)$ is normal in $\mf G$, it is normalized by $\mf H$. Hence 
  the group generated by $\mf H_{S}$ and $\mc R_u(\mf G)_{S}$ is amenable.  
\end{proof}

From now on, $\mf G$ will be reductive and $\mf H \leqslant \mf G$ be $S$-maximal amenable. Consider the Levi decomposition $\mf H = \mf H^R \ltimes \mc R_u(\mf H)$.

\begin{Lemma}\label{MA absorbs unipotent}
 Assume that $\mf H\leqslant \mf G$ is $S$-maximal amenable. Let $\mf U\leqslant \mf G$ be a unipotent subgroup and assume that $\mf H\leqslant \mc N(\mf U)$. Then $\mf U\leqslant\mc R_u(\mf H)\leqslant\mf H$. In particular, given any subgroup $\mf P\leqslant \mf G$, if $\mf H\leqslant\mf P\leqslant\mf G$, then $\mc R_u(\mf H)\geqslant \mc R_u(\mf P)$.
\end{Lemma}
\begin{proof}
  Assume that $\mf H$ is $S$-MA and denote by $\mf H'$ the group generated by $\mf H$ and $\mf U$. Observe that $\mf U\leqslant\mf H'$ is normal. The group $\mf H'_{S}/\mf U_{S}=\mf H_{S}/(\mf U_{S}\cap \mf H_{S})$ is amenable and hence $\mf H'_{S}$ is. Therefore $\mf H'$ is $S$-amenable and the maximality implies that $\mf H'=\mf H$. The subgroup $\mf U\leqslant \mf H$ is a normal unipotent subgroup which has to be contained in the unipotent radical.
\end{proof}

\begin{Prop}\label{prop free parabolic}
  Assume that $\mf H\leqslant \mf G$ is $S$-maximal amenable. Then 
  \[
  \mc R_u(\mc N_{\mf G}(\mc R_u(\mf H)))=\mc R_u(\mf H),
  \]
 and moreover $\mf P\coloneqq\mc N_{\mf G}(\mc R_u(\mf H))$ is a parabolic subgroup of $\mf G$ containing $\mf H$. 
\end{Prop}
\begin{proof}
  Set $\mf U\coloneqq\mc R_u(\mf P)$. By definition, $\mc R_u(\mf H)$ is a unipotent group normalized by $\mf P$, so it is contained in its unipotent radical, $\mc R_u(\mf H)\leqslant \mf U$. 
  Since $\mc R_u(\mf H)$ is normal in $\mf H$, we have that $\mf H\leqslant \mc N_{\mf G}(\mc R_u(\mf H))=\mf P$. Lemma \ref{MA absorbs unipotent} yields that $\mc R_u(\mf H)\geqslant \mc R_u(\mf P)=\mf U = \mc R_u(\mc N_{\mf G}(\mc R_u(\mf H)))$. 

  The moreover part follows from \cite[Corollary 30.3.B]{HumphreysLinear1975}.
\end{proof}

\begin{Lemma}[{\cite[Proposition I.1.3.2, Proposition I.1.6.3, Proposition I.2.3.6]{MargulisDiscrete1991}}]\label{lem:compact=amenable=anisotropic}
Let $k$ be a local field. If $\mf G$ is a semisimple $k$-group, then the following conditions are equivalent:
\begin{enumerate}
    \item $\rk_k \mf G = 0$;
    \item $\mf G(k)$ is compact;
    \item $\mf G(k)$ doesn't contain any unipotent elements;
    \item $\mf G(k)$ is amenable.
\end{enumerate}
If $\mf G$ is a reductive $k$-group, then the first two conditions are equivalent and the last two conditions are equivalent. 
\end{Lemma}

Note that the derived subgroup of a connected reductive group is always semisimple \cite[Proposition 2.2]{BorelGroupes1965}. Therefore Lemma \ref{lem:compact=amenable=anisotropic} implies the following.

\begin{Lemma}\label{lemma:no-unipotents}
   Let $\mf H$ be a connected $K$-group with Levi decomposition $\mf H\cong \mf H^R\ltimes \mc R_u(\mf H)$. Then the following conditions are equivalent:
  \begin{enumerate}
      \item $\mf H$ is $S$-amenable;
      \item $\rk_{K_v} \mf H^R = \rk_{K_v} \mc Z(\mf H^R)$ for all $v\in S$;
      \item $\mc D (\mf H^R)$ is $S$-compact.
  \end{enumerate}
In particular, if $\mf H$ is $S$-amenable, then all unipotent elements of $\mf H_{S}$ are contained in $\mc R_u(\mf H)_{S}$.
\end{Lemma}

For infinite $S$, we can say much more.

\begin{Lemma}\label{lemma:no-unipotents infinite S}
   Let us assume that $S$ is infinite.
   Let $\mf H$ be a connected $K$-group with Levi decomposition $\mf H\cong \mf H^R\ltimes \mc R_u(\mf H)$. Then $\mf H$ is $S$-amenable if and only if $\mf H^R$ is virtually abelian.  
\end{Lemma}
\begin{proof}
  The lemma follows from Lemma \ref{lemma:no-unipotents} once we observe that the set of places $v\in S$ for which $\mc D(\mf H^R)(K_v)$ is compact is finite \cite[I.3.2.3]{MargulisDiscrete1991}.
\end{proof}

Therefore, we can concentrate on understanding maximal amenable subgroups without unipotent elements, that is, maximal amenable reductive subgroups $\mf H \leqslant \mf G$ in a reductive group $\mf G$. 

\begin{Lemma}\label{lem:ranks-and-central-tori}
Let $\mf H$ be a reductive $S$-maximal amenable subgroup of a reductive group $\mf G$ and $\mf Z$ be the maximal $K$-split torus of $\mc Z(\mf G)$. Then the following holds:
\begin{enumerate}
\item \label{item: rank central kv torus} 
for every $v\in S$, let $\mf T_v$ be a $K_v$-defined maximal $K_v$-split torus of $\mf H$. Then $\mf T_v$ is central in $\mf H^\circ$.
\item \label{item: k-split torus of H is central}
if $\mf H$ is connected, then $\rk_K \mf H = \rk_K \mf Z$. In particular, in this case $\mf Z$ is the maximal $K$-split torus of $\mf H$.
\end{enumerate}
\end{Lemma}
\begin{proof}
For every $v\in S$, Lemma \ref{lemma:no-unipotents} yields that $\mf H^\circ(K_v)$ is the almost direct product of its center with a compact subgroup. Therefore any $K_v$-defined $K_v$-split subgroup of $\mf H^\circ(K_v)$ is contained in the center, so  \ref{item: rank central kv torus} is proven. 

Clearly $\mf H\geqslant \mf Z$. Assume that we have a $K$-split torus $\mf Z'>\mf Z$ properly containing $\mf Z$. By the first point, $\mf Z'$ is in the center of $\mf H$. On the other hand $\mf Z'$ is not central in $\mf G$, therefore it has at least one non-trivial root. Then by \cite[Proposition I.1.1.3]{MargulisDiscrete1991}, \cite[Corollaire 3.18]{BorelGroupes1965}, $\mf Z'$ normalises some $K$-unipotent subgroup $\mf U$ in $\mf G$. Since $\mf Z'$ is central in $\mf H$, we have $\mf H\leqslant \mc C_{\mf G}(\mf Z')$. This centralizer is contained in a parabolic subgroup $\mf P$ with non-trivial unipotent radical $\mf U$. But in this case $\mf H$, being maximal amenable, has to contain $\mf U$ by Lemma \ref{MA absorbs unipotent}. This contradicts the assumption that $\mf H$ is reductive and thus proves \ref{item: k-split torus of H is central}.
\end{proof}

Note that Item \ref{item: k-split torus of H is central} is not true if $\mf H$ is not connected. Indeed the normalizer of the $\mb Q$-split torus of $\mf{SL}_3(\mb R)$ is maximal amenable and does not satisfy Item \ref{item: k-split torus of H is central}, see Example \ref{ex: normalizer split torus is MA} below.

\begin{Lemma}\label{lem:maximal-amenable-in-reductive}
Let $\mf G$ be a reductive group and $\mf H \leqslant \mf G$ be a reductive $S$-maximal amenable subgroup. Let $\mf T$ be the central torus of $\mf H^\circ$. Then $\mf H \leqslant\mc N_{\mf G}(\mf T)$. In particular,
\begin{itemize}
    \item $\mf H/\mf T \leqslant \mc N_{\mf G}(\mf T)/\mf T$ is a semisimple $S$-maximal amenable subgroup of the semisimple group $\mc N_{\mf G}(\mf T)/\mf T$;
    \item for every $v\in S$, $\rk_{K_v}(\mf H/\mf T)=0$;
    \item if $S$ is finite, the group $(\mf H/\mf T)_S$ is compact;
    \item if $S$ is infinite, $(\mf H/\mf T)_S$ is finite.
\end{itemize}
\end{Lemma}
\begin{proof}
Let $\mf H\leqslant\mf G$ be reductive and $S$-maximal amenable and $\mf T \leqslant \mf H$ be the central torus of $\mf H^\circ$. Then $\mf H^\circ \leqslant \mc C_{\mf G}(\mf T)$. As $\mf H^\circ$ is normal in $\mf H$ and every automorphism of $\mf H^\circ$ has to preserve the central torus, we have that $\mf H\leqslant\mc N_{\mf G}(\mf T)$ which is a finite index overgroup of the Levi subgroup $\mc C_{\mf G}(\mf T)$. Moreover, $\mf H$ normalizes $\mc C_{\mf G}(\mf T)$ and therefore its center $\mf T'\geqslant \mf T$. Since $\mf H$ is $S$-maximal amenable, it contains $\mf T'$ which is therefore central in $\mf H^\circ$; hence $\mf T' = \mf T$ which ensures that $\mc N_{\mf G}(\mf T)/\mf T$ is semisimple. The rest follows from Lemmas \ref{lem:compact=amenable=anisotropic}, \ref{lemma:no-unipotents} and \ref{lemma:no-unipotents infinite S}.
\end{proof}

\begin{Rem}
The connected component of the quotient $\mf H/\mf T$ can still be viewed as a subgroup of $\mf G$ since by the structure theory of reductive groups $\mf H$ is an almost direct product of its central torus $\mf T$ and $\mc D(\mf H)$.
\end{Rem}

The above results allow us to deduce maximal amenability in certain important cases.

\begin{Def}
An inclusion $\mf H \leqslant \mf G$ is $S$-cocompact if for every $v\in S$ the group $\mf H(K_v)$ is cocompact in $\mf G(K_v)$. 
\end{Def}

\begin{Def}\label{def:s-ample}
A $K$-torus $\mf T$ in $\mf G$ is said to be $S$-\textbf{ample} if the following conditions hold:
\begin{enumerate}
    \item \label{item: s-ample rank center} $\rk_K(\mf T)=\rk_K(\mc Z(\mf G))$;
    \item\label{item: s-ample rank centralizer} $\mf T$ is $S$-cocompact in $\mc N_{\mf G}(\mf T)$;
    \item\label{item: s-ample no split compact} $\mf T$ has no non-trivial $S$-cocompact subgroups.
\end{enumerate}
\end{Def}

By Lemma \ref{lem:compact=amenable=anisotropic}, Item \ref{item: s-ample rank centralizer} above is equivalent to the following. 

\begin{enumerate}
    \item[\mylabel{item: s-ample rank centralizer'}{(ii')}] for every $v\in S$ we have $\rk_{K_v}\mf T = \rk_{K_v}\mc C_{\mf G}(\mf T)$.
\end{enumerate}

Similarly, Item \ref{item: s-ample no split compact} is equivalent to

\begin{itemize}
    \item[\mylabel{item: s-ample no split compact'}{(iii')}] for every proper $K$-defined sub-torus $\mf T'\leqslant \mf T$ there exists $v\in S$ such that $\rk_{K_v}(\mf T')< \rk_{K_v}(\mf T)$.
\end{itemize}     
        
Let us also remark that if $\mf T$ is an $S$-ample torus, then $\mf T$ is the smallest $K$-defined torus in $\mf G$ containing the maximal $K_v$-split subtori of $\mf T$ for all $v\in S$.

\begin{Lemma}\label{lem:s-ample-unique}
  Let $\mf G$ be a reductive group and let $\mf T\leqslant \mf G$ be an $S$-ample torus. Then $\mf T$ is the unique $S$-ample sub-torus of $\mc C_{\mf G}(\mf T)$. 
\end{Lemma}
\begin{proof}
Let $\mf T_0\leqslant \mc C_{\mf G}(\mf T)$ be an $S$-ample torus. Denote by $\mf H\leqslant \mc C_{\mf G}(\mf T)$ be an almost complement of $\mf T$, that is a subgroup such that $\mf H$ and $\mf T$ generate $\mc C_{\mf G}(\mf T)$ and have finite intersection. By \ref{item: s-ample rank centralizer'}, \[\rk_{K_v}(\mf H\cap \mf T_0)\leqslant \rk_{K_v}(\mf H)=0\] for every $v\in S$. This implies that $\rk_{K_v}(\mf T_0\cap \mf T)=\rk_{K_v}(\mf T_0)$ for every $v\in S$. By \ref{item: s-ample no split compact'}, we deduce that $\mf T_0\subseteq \mf T$. Now by \ref{item: s-ample rank centralizer'} applied to $\mf T_0$ we deduce that $\rk_{K_v} \mf T = \rk_{K_v}\mf T_0$ for every $v\in S$, and then by \ref{item: s-ample no split compact'} it follows that $\mf T= \mf T_0$.
\end{proof}

\begin{Prop}\label{prop:normalizer-max-amenable-sufficient}
Let $\mf G$ be a reductive group and $\mf T \leqslant \mf G$ be an $S$-ample torus. Let $\mf H$ be an $S$-amenable subgroup containing $\mf T$. Then $\mf H \leqslant \mc N_{\mf G}(\mf T)$.
\end{Prop}
\begin{proof}
Since $\mf T$ is connected, $\mf T\leqslant \mf H^\circ$. 
Let us first prove that $\mf H$ cannot contain unipotent elements. Indeed, if it does, then the unipotent radical $\mc R_u(\mf H)$ is nontrivial (Lemma \ref{lemma:no-unipotents} and Lemma \ref{lemma:no-unipotents infinite S}) and is normalised by $\mf T$. In particular, $\mf T$ is contained in a proper parabolic subgroup of $\mf G$ and therefore $\mc C_{\mf G}(\mf T)$ has to contain a $K$-split subtorus $\mf T_0$ which is not central in $\mf G$.
By condition \ref{item: s-ample rank centralizer} in Definition \ref{def:s-ample}, this torus $\mf T_0$ is contained in $\mf T$; this is, however, in contradiction with the condition \ref{item: s-ample rank center} in Definition \ref{def:s-ample}. 

Thus, $\mf H$ is reductive and Lemma \ref{lem:maximal-amenable-in-reductive} applies. Let $\mf T'$ be the central torus of $\mf H^\circ$; by construction, $\mf T'$ commutes with $\mf T$. By Lemma \ref{lem:maximal-amenable-in-reductive}, for every $v\in S$
\[
\rk_{K_v}(\mf T/(\mf T\cap \mf T')) \leqslant \rk_{K_v}(\mf H/\mf T')= 0.
\] 
However, condition \ref{item: s-ample no split compact} in Definition \ref{def:s-ample} implies that $\mf T$ has no non-trivial $S$-anisotropic quotients, hence $\mf T \leqslant \mf T'$, which implies $\mc C_{\mf G}(\mf T)\geqslant \mc C_{\mf G}(\mf T') \geqslant \mf H^\circ$. By Lemma \ref{lem:s-ample-unique} the torus $\mf T\leqslant \mf H^\circ$ is the unique $S$-ample torus and hence it is normalized by $\mf H$. That is, $\mf H\leqslant \mc N_{\mf G}(\mf T)$, and we are done.
\end{proof}

The discrepancy between the necessary condition of Lemma \ref{lem:maximal-amenable-in-reductive} and the sufficient condition of Proposition \ref{prop:normalizer-max-amenable-sufficient} is the fact that -- to the best of authors' knowledge -- maximal $K$-defined $S$-compact subgroups of reductive groups are not well-understood. However, if we allowed to work up to ``$S$-cocompact inclusions'', the situation becomes easier. 

\begin{Def}
Two subgroups $\mf H$ and $\mf K$ of $\mf G$ are said to be $S$-\textbf{cocompactly equivalent} if their intersection is $S$-cocompact in both of them.
\end{Def}

Clearly all groups in the same $S$-cocompact equivalence class are either $S$-amenable or not. Thus, one can speak of amenable and non-amenable $S$-cocompact equivalence classes.

\begin{Def}\label{def: S cocompaclty maximal amenable}
An $S$-maximal amenable subgroup $\mf H\leqslant \mf G$ is $S$-\textbf{cocompactly maximal amenable} if its $S$-cocompact equivalence class is maximal.
\end{Def}

Let us first observe that Proposition \ref{prop:normalizer-max-amenable-sufficient} implies the following. 

\begin{Cor}\label{cor: S-ample is S CMA}
 The normalizer of an $S$-ample torus $\mf T\leqslant \mf G$ in a reductive group $\mf G$ is $S$-cocompactly maximal amenable. 
\end{Cor}

\begin{Thm}\label{thm:maximal-amenable-algebraic-subgroups}
Let $\mf G$ be a connected $K$-group. Let $S\subseteq \ms R$ be a set of valuations. 
\begin{enumerate}
\item\label{mainthm: S maximal amenable} Let $\mf H\leqslant \mf G$ be an $S$-maximal amenable subgroup with Levi decomposition $\mf H=\mf H^R\ltimes \mf U$, where $\mf H^R$ is reductive and $\mf U\coloneqq \mc R_u(\mf H)$. 
Then $\mf U$ is the unipotent radical of a $K$-parabolic subgroup $\mf P\coloneqq \mc N_{\mf G}(\mf U)$  of $\mf G$. The parabolic subgroup $\mf P$ has a Levi decomposition $\mf P = \mf L\ltimes \mf U$ such that $\mf H^R \leqslant \mf L$. The central torus $\mf T$ of $\mf H^R$ contains $\mc Z(\mf L)^\circ$ and the image of $\mf H^R$ in $\mc N_{\mf L}(\mf T)/\mf T$ is maximal among semisimple $S$-compact $K$-defined subgroups.

Moreover, the following properties hold:
\begin{enumerate}
  \item\label{item1: unipotent elements in radical} $\mf U$ is a maximal unipotent group normalized by $\mf H^R$ and all unipotent elements in $\mf H_S$ belong to $\mf U_S$.
  \item\label{item1: almost direct product} $(\mf H^R)_S$ is an almost direct product of an abelian group and a compact group; if $S$ is infinite, it is virtually abelian.  
  \item \label{item1: rank central kv torus} For every $v\in S$, let $\mf T_v$ be a $K_v$-defined maximal $K_v$-split torus of $\mf H$. Then $\mf T_v$ is central in $\mf H^\circ$; in particular, $\mf T_v$ is unique.
\item \label{item1: k-split torus of H is central}
If $\mf H$ is connected, then $\rk_K \mf H = \rk_K \mc Z(\mf L)^\circ$. In particular, the maximal $K$-split torus of $\mc Z(\mf L)^\circ$ is the maximal $K$-split torus of $\mf H$.
\end{enumerate}

\item \label{mainthm: S cocompactly maximal amenable}
Every $S$-cocompactly maximal amenable subgroup $\mf H\leqslant \mf G$ is $S$-cocompactly equivalent to a subgroup $\mf H_{\mf T}\leqslant \mf H$ of the form
\[
\mf H_{\mf T} \coloneqq \mf T \ltimes \mf U,
\]
where $\mf U$ is the unipotent radical of the $K$-parabolic subgroup $\mf P\coloneqq \mc N_{\mf G}(\mf U)$ of $\mf G$ with a Levi decomposition $\mf P = \mf L\ltimes \mf U$, and $\mf T\leqslant \mf L$ is an $S$-ample torus. 

Moreover if we denote by $\mf T_a$ the maximal $K$-anisotropic subtorus of $\mf T$, we have
\begin{enumerate}
  \item\label{item2: unique smallest parabolic} $\mc N_{\mf G}(\mf U)$ is the unique smallest parabolic subgroup of $\mf G$ containing $\mf T_a\ltimes \mf U$.
  \item\label{item2: normalizer is MA} $\mc N_{\mf G}(\mf H_{\mf T})=\mc N_{\mf L}(\mf T)\ltimes \mf U$ is the unique $S$-maximal amenable subgroup of $\mf G$ containing $\mf T_a\ltimes \mf U$.
  \item\label{item2: T and U uniquely defined} The subgroup $\mf U$ is uniquely determined by the equivalence class of $\mf H_{\mf T}$ and the subgroup $\mf T$ is uniquely determined up to conjugation by $\mf U$.
  \item\label{item2: minimal and maximal class} $\mf H_{\mf T}$ is the minimal $K$-subgroup in its equivalence class and $\mc N_{\mf G}(\mf H_{\mf T})$ is the maximal $K$-subgroup in the equivalence class.
\end{enumerate}
In particular, for every subgroup $\mf H_{\mf T}$ as above the normalizer $\mc N_{\mf G}(\mf H_{\mf T})$ is $S$-cocompactly maximal amenable.
\end{enumerate}
\end{Thm}
\begin{proof}
The main claim of part \ref{mainthm: S maximal amenable} is a direct consequence of the reduction steps contained in Lemma \ref{lem:ma-contains-unipotent-radical} till Lemma \ref{lem:maximal-amenable-in-reductive}. Indeed, by Lemma \ref{lem:ma-contains-unipotent-radical} $\mf H\geqslant\mc R_u(\mf G)$ and hence by passing to the quotient, we can assume that $\mf G$ is reductive; then Proposition \ref{prop free parabolic} yields that $\mc N_{\mf G}(\mf U)$ is a parabolic group and by Lemma \ref{lemma:no-unipotents}, all unipotent elements of $\mf H_S$ are contained in $\mf U$. Finally $\mf H^R$ is contained in a Levi complement $\mf L$ of $\mf U$ in $\mf P$ and it is an $S$-maximal amenable reductive group; we can conclude by Lemma \ref{lem:maximal-amenable-in-reductive}.
Item \ref{item1: unipotent elements in radical} and Item \ref{item1: almost direct product} follow from Lemma \ref{lemma:no-unipotents}, and the remaining two points restate Lemma \ref{lem:ranks-and-central-tori}.
  
Let us now prove part \ref{mainthm: S cocompactly maximal amenable}. 
Part \ref{mainthm: S maximal amenable} gives us information about the structure of $\mf H$; let $\mf T'$ be the central torus of the respective subgroup $\mf H^R\leqslant \mf L$ and $\mf U\leqslant \mf P$ be the respective unipotent subgroup. Let $\mf T\leqslant \mf T'$ be the minimal $S$-cocompact subtorus of $\mf T'$. Then $\mf H_{\mf T}=\mf T \ltimes \mc R_u(\mf H)$ is $S$-cocompactly equivalent to $\mf H$. 
Let us prove that the torus $\mf T\leqslant \mf L$ is $S$-ample. By Item \ref{item1: k-split torus of H is central}, $\rk_K\mf T = \rk_K \mc Z(\mf L)^\circ$ and $\mc C_{\mf G}(\mf T) = \mc C_{\mf L}(\mf T)$; by construction of $\mf T$ it doesn't have nontrivial $S$-cocompact subgroups. Finally, if $\rk_{K_v}\mc C_{\mf L}(\mf T) > \rk_{K_v}\mf T$ for some particular $v\in S$, then by \cite[Corollary 3.2]{PrasadWeakly2009} there is a $K$-defined subtorus $\mf T''\geqslant \mf T$ in $\mc C_{\mf L}(\mf T)$ such that $\rk_{K_v} \mf T'' > \rk_{K_v}\mf T$. Since $\mf T''$ still normalizes $\mc R_u(\mf H)$, the subgroup $\mf T''\rtimes \mc R_{u}(\mf H)$ is amenable and $\mf T\rtimes \mc R_{u}(\mf H)$ is not $S$-cocompact in it, contradicting the assumption that $\mf H$ is $S$-cocompactly maximal amenable. Therefore $\mf T$ is $S$-ample in $\mf L$.

Since all $S$-cocompactly maximal amenable subgroups of $\mf G$ contain $\mc R_u(\mf G)$, for the rest of the proof we assume that $\mf G$ is reductive.

Let us show that Item \ref{item2: unique smallest parabolic} holds. Let $\mf P'\geqslant\mf T_a\ltimes\mc R_u(\mf H)$ be a parabolic subgroup. Since $\mf P''\coloneqq \mf P'\cap \mf P\geqslant\mf U$, \cite[Proposition 4.4. b)]{BorelGroupes1965} implies that this intersection is a parabolic subgroup of $\mf G$. 
By \cite[Proposition 4.4. c)]{BorelGroupes1965} the group $\mf L\cap \mf P''$ is a parabolic subgroup of $\mf L$ containing $\mf T_a$. Remark now that $\mf T_a$ is $S$-ample in $\mc D(\mf L)$, so Corollary \ref{cor: S-ample is S CMA} implies that $\mf T_a$ is $S$-cocompactly maximal amenable in $\mc D(\mf L)$. Therefore $\mf T_a$ cannot normalize any unipotent subgroup of $\mc D(\mf L)$ and hence of $\mf L$. That is, $\mc R_u(\mf P'')=\mf U$ so that $\mf P''=\mf P$. So $\mf P'\geqslant \mf P$ therefore the proof of Item \ref{item2: unique smallest parabolic} is completed.

Let us prove Item \ref{item2: normalizer is MA}.
Let $\mf H'$ be an $S$-maximal amenable subgroup containing $\mf T_a\ltimes \mf U$. By Part \ref{mainthm: S maximal amenable} there is a torus $\mf T'$ and a unipotent group $\mf U'$ such that $\mf T'\ltimes\mf U'\leqslant\mf H'$ is $S$-cocompact and $\mf U'\geqslant\mf U$. We also know that $\mf P'\coloneqq \mc N_{\mf G}(\mf U')$ is a parabolic subgroup containing $\mf H'$. By Item  \ref{item2: unique smallest parabolic} we have that $\mf P'\geqslant \mf P$. Since $\mf U'\geqslant\mf U$, \cite[Proposition 4.4. c)]{BorelGroupes1965} implies that $\mf P' = \mf P$. Observe that $\mf H'\geqslant \mc Z(\mf L)^\circ$ and hence $\mf H'\geqslant \mf T$. Now by passing to the quotient by $\mf U$, we can apply Proposition \ref{prop:normalizer-max-amenable-sufficient} to deduce that $\mf H'\leqslant \mc N_{\mf L}(\mf T)\ltimes \mf U$. Finally since $\mf T$ is $S$-ample, we get that $\mf H_{\mf T}$ is $S$-cocompact in $\mc N_{\mf L}(\mf T)\ltimes \mf U$, so it is $S$-cocompactly maximal amenable as claimed and $\mc N_{\mf L}(\mf T)\ltimes \mf U$ is the unique $S$-maximal amenable group containing $\mf T_a\ltimes \mf U$.

In Item \ref{item2: T and U uniquely defined}, the uniqueness of $\mf U$ is clear (it is the unipotent radical of every subgroup in the $S$-cocompact equivalence class), and the uniqueness of $\mf T$ follows form Lemma \ref{lem:s-ample-unique}. Finally Item \ref{item2: minimal and maximal class} is just a consequence of the previous items. 

If $\mf T$ is $S$-ample, the subgroup $\mf H_{\mf T}$ doesn't have any proper $S$-cocompact subgroups. Thus, the fact that the normalizer of every subgroup $\mf H_{\mf T}$ as above is $S$-cocompactly maximal amenable follows from Item \ref{item2: normalizer is MA}.
\end{proof}

Let us now show that every $S$-amenable group is $S$-cocompactly contained in an $S$-cocompactly maximal amenable group.

\begin{Prop}\label{prop: every sma is contained scma}
  Let $\mf H\leqslant\mf G$ be an $S$-maximal amenable subgroup. Then there is an $S$-cocompactly maximal amenable subgroup $\mf H' = \mc N_{\mf G}(\mf T')\ltimes \mc R_u(\mf H) \leqslant \mf G$ for some torus $\mf T'\leqslant \mf G$ such that 
  \begin{itemize}
      \item $\mf H\cap \mf H'$ is $S$-cocompact in $\mf H$;
      \item $\mf H'/(\mf H\cap \mf H')=\mf T'/(\mf H\cap \mf T')$ is a torus of with $K$-rank $0$;
      \item $\mc R_u(\mf H)=\mc R_u(\mf H')$. 
  \end{itemize}
\end{Prop}
\begin{proof}
  By Theorem \ref{thm:maximal-amenable-algebraic-subgroups} \ref{mainthm: S maximal amenable}, we have that $\mf H=\mf H^R\ltimes \mc R_u(\mf H)$ and there is a torus $\mf T\leqslant \mf H^R$ which is $S$-cocompact so that $\mf T\ltimes \mc R_u(\mf H)$ is $S$-cocompact in $\mf H$. We can furthermore assume that $\mf T$ has no non-trivial $S$-cocompact subgroups. Denote by $\mf Z\leqslant \mf T$ the maximal $K$-split subtorus. Then $\mc C_{\mf G}(\mf T)$ is a reductive group. If $S$ is finite, by \cite[Corollary 3.2]{PrasadWeakly2009}, there is a torus $\mf T'\geqslant \mf T$ such that $\rk_K(\mf T')=\rk_K(\mf T)$ and such that $\rk_{K_v}(\mc C_{\mf G}(\mf T))=\rk_{K_v}(\mf T')$ for every $v\in S$. Clearly $\mf T'$ contains an $S$-ample torus of $\mc C_{\mf G}(\mf Z)$, therefore Theorem \ref{thm:maximal-amenable-algebraic-subgroups} \ref{mainthm: S cocompactly maximal amenable} implies that $\mc N_{\mf G}(\mf T')\ltimes \mf U$ is $S$-cocompactly maximal amenable. 
  
  Finally if $S$ is infinite, then we can iterate the above construction with the finite subsets of $S$. These iterations have to stop after finitely many steps and gives us the desired result. 
\end{proof}

Let us remark that there is a natural class of examples of $S$-maximal amenable subgroups which are not $S$-cocompactly maximal amenable.

\begin{Ex}\label{ex: normalizer split torus is MA}
Let $\mf G$ be reductive $K$-split group (that is $\rk_K \mf G =\rk_{\bar K} \mf G$), $\mf S$ a maximal $K$-split torus in $\mf G$. Then $\mc N_{\mf G}(\mf S)$ is $S$-maximal amenable and does not contain any unipotents; in particular it is not cocompactly maximal amenable.

Indeed, if $\mc N_{\mf G}(\mf S)$ normalizes a unipotent subgroup $\mf U$, then $\mf U$ is contained in a minimal parabolic $\mf P$ which has $\mf S$ as its central torus. However, then there exists an element $n\in \mc N_{\mf G}(\mf S)(K)$ which conjugates $\mc R_u(\mf P)$ to $\mc R_u(\mf P^-)$, where $\mf P^-$ is the opposite parabolic. However, $\mc R_u(\mf P)\cap \mc R_u(\mf P^-) = \{1\}$, yielding $\mf U$ is trivial. An application of Lemma \ref{lem:maximal-amenable-in-reductive} concludes.
\end{Ex}

Let us remark that all results of this section remain valid if $K$ is assumed to be a local field with $S$ consisting of its absolute value. In this case, $S$-amenability amounts to amenability with respect to the topology coming from this absolute value.

In the case $K=\mathbb R$, the first part of Theorem \ref{mainthm: S cocompactly maximal amenable} then amounts to the classification of maximal amenable subgroups obtained by Moore in \cite{MooreAmenable1979}. Indeed, in this case the Zariski closure of an amenable subgroup of $\mf G(\mb R)$ is amenable \cite[Theorem 2.11]{MooreAmenable1979}. Note that the example above isn't included into Moore's classification; indeed, it fails the condition to be isotropically connected from \cite{MooreAmenable1979}. For other local fields our technique allows to deduce the description previously observed in \cite[Corollary 1.18]{BaderAlmost2017}. 

\begin{Prop}
  Let $K$ be a local field and $\mf G$ a reductive $K$-group. Let $H\leqslant \mf G(K)$ be a maximal amenable subgroup with respect to the topology of $K$ such that its Zariski closure $\mf H$ is (Zariski) connected. Then there is a $K$-defined parabolic subgroup $\mf P\leqslant \mf G$ containing $H$ and such that $\mc R(\mf P)\leqslant H$ is cocompact. If moreover $H$ is cocompactly maximal amenable, then $H=\mf H(K)$. 
\end{Prop}
\begin{proof}[Sketch of proof]
In this proof we will tacitly identify all algebraic subgroups with their $K$-points.

Let $\mf H$ be the Zariski closure of $H$; note that $\mf H$ is not necessarily amenable. Since $H$ is maximal amenable, $\mc R_u(\mf H)\leqslant H$. By construction, $H\leqslant \mc N(\mc R_u(\mf H))\eqqcolon \mf P$. Arguing as in Proposition \ref{prop free parabolic}, we deduce that $\mf P$ is a parabolic subgroup. Clearly, $\mc R(\mf P)\leqslant \mf H$. We are left to show that the quotient $H/\mc R(\mf P)$ is compact. It is contained in the reductive group $\mf H/\mc R(\mf P)$ and it decomposes as an almost direct product $\mf C \times \prod_{i=1}^k \mf H_i$ where $\mf C$ is the center of $\mf H$ and $\mf H_i$ are almost simple. Since $H$ is maximal amenable, $H/\mc R(\mf P)$ decomposes as an almost direct product $\mf C\times \prod_{i=1}^k H_i$ with $H_i\leqslant \mf H_i$ Zariski dense. By \cite[Corollary 3.2.19]{ZimmerErgodic1984}, we see that $H_i$ is compact. If $H$ is furthermore cocompactly maximal amenable, then $\mf H_i$ has to be compact and hence $\mf H_i=H_i$.
Since $\mf H$ is connected, we can argue as in Lemma \ref{lem:ranks-and-central-tori} \ref{item1: k-split torus of H is central}, to deduce that $\mf C$ is compact. This finishes the proof.
\end{proof}

\section{$K(S)$-points}
In this section, $K$ is a global field, $\mf G \leqslant \mf{GL}_n$ is a connected algebraic group. We will set $\mf G(K(S)):=\mf G\cap \mf{GL}_n(K(S))$. 

We will denote by $\mathrm X_K(\mf G)$ the set of $K$-characters of the group, that is the set of maps of algebraic groups from $\mf G$ to $K$ seen as the additive $K$-group. We let $\ms T(\mf G)\subset \ms R$ be the subset of valuations for which $\mf G$ is anisotropic, that is those $v\in\ms R$ such that $\mf G(K_v)$ is compact. 

By \cite[I.3.2.3, Theorem I.3.2.4 (a), I.3.2.6]{MargulisDiscrete1991} if $\mf G$ is a connected algebraic $K$-group and $\ms R_\infty\setminus \ms T(\mf G)\subseteq S\subseteq \ms R$, then the subgroup $\mf G(K(S))$ is discrete in $\mf G_S$, and we will assume it for the rest of the section. Moreover, $\mf G(K(S))$ is a lattice in $\mf G_S$ if and only if $\mathrm X_K(\mf G)=1$.


\begin{Lemma}\label{lemma:cocompact implies Zariski dense}
 Let $\mf G$ be an algebraic group without $S$-cocompact subgroups. Assume that $\Gamma\leqslant\mf G(K)$ is an $S$-cocompact subgroup. Then $\Gamma$ is Zariski dense in $\mf G$.  
\end{Lemma}
\begin{proof}
    The Zariski closure of $\Gamma$ is also $S$-cocompact in $\mf G$ and $K$-defined. 
\end{proof}

\begin{Rem}\label{rem:HT(K(S)) Zariski dense in HT}
In the notation of Theorem \ref{thm:maximal-amenable-algebraic-subgroups} \ref{mainthm: S cocompactly maximal amenable}, it follows from the lemma that if $\mf H$ is an $S$-commensurably maximal amenable subgroup, then the Zariski closure of any finite index subgroup of $\mf H_{\mf T}(K(S))$ contains $\mf T_a\ltimes \mf U$. Indeed, $(\mf T_a\ltimes \mf U)(K(S))$ (and any finite index subgroup of it) is $S$-cocompact in $(\mf T_a\ltimes \mf U)_S$ by \cite[Theorem I.3.2.4, Remark I.3.2.6]{MargulisDiscrete1991} and $\mf T_a\ltimes \mf U$ has no $K$-defined $S$-cocompact subgroups. 
\end{Rem}

\begin{Lemma}\label{lem: amenable iff K(S) points amenable}
  Let $\mf H\leqslant \mf G$ be a subgroup. Then $\mf H$ is $S$-amenable if and only if $\mf H(K(S))$ is amenable.
\end{Lemma}
\begin{proof}
  If $\mf H_S$ is amenable, then its closed subgroup $\mf H(K(S))$ will be amenable. Conversely assume that $\mf H(K(S))$ is amenable. Since amenability is preserved in the commensurability class, we can assume that $\mf H$ is connected. Consider the Levi decomposition $\mf H\cong \mf H^R\ltimes \mc R_u(\mf H)$. Then $\mf H^R(K(S))<\mf H(K(S))$ is amenable, and so therefore is $\mc D (\mf H^R)(K(S))$  By \cite[Theorem I.3.2.5]{MargulisDiscrete1991}, this implies that $\mf H^R_S$ is amenable.
\end{proof}

\begin{Prop}[Poor man's Tits alternative]\label{prop: poor mans tits alternative}
  Let $K$ be a global field and $\Lambda\leqslant \mf G(K)$ be an amenable subgroup. Then the Zariski closure of $\Lambda$ is virtually solvable.  
\end{Prop}
\begin{proof}
 Let $\mf H$ be the Zariski closure of $\Lambda$. We want to prove that the semisimple group $\mf H' \coloneqq \mf H/\mc R(\mf H)$ is finite. Since $\mf H'$ is an almost direct product of almost simple algebraic groups, replacing $\Lambda/(\Lambda \cap \mc R(\mf H))$ by a finite extension of it, we can assume that $\mf H = \mf H_1\times\dots\times\mf H_r$. We have to prove the finiteness of every $\mf H_i$. 

 We let $\Lambda_i\leqslant\mf H_i(K)$ be the image of $\Lambda$ in $\mf H_i$; it is an amenable Zariski dense subgroup. Let $\mf P_i<\mf H_i$ be a minimal $K$-parabolic subgroup. Then the quotient $\mf H_i(K_v)/\mf P_i(K_v)$ is compact for every $v\in S$. The group $\Lambda_i\leqslant \mf H_i(K_v)$ is amenable and Zariski dense, so \cite[Corollary 3.2.19]{ZimmerErgodic1984} implies that $\Lambda_i\leqslant \mf H_i(K_v)$ is precompact. Then we deduce that $\Lambda_i$ is contained in a restricted product of compact subgroups $\prod_{v\in \ms R}M_v \leqslant \mf H_i(\mb A_K)$; this restricted product is a direct limit of an increasing sequence of compact groups. Since $\Lambda_i < \mf H_i(K)$ and the latter is discrete in $\prod_{v\in \ms R}\mf H_i(K_v)$, we deduce that $\Lambda_i$ is locally finite. Now, the Jordan--Schur theorem implies \cite[Corollary 12.1.12]{PassmanAlgebraic1977} that $\Lambda_i$ is virtually abelian. It follows that $\Lambda_i$ is finite, since its Zariski closure $\mf H_i$ is supposed to be semisimple and thus its identity component can't be abelian unless it is trivial.
\end{proof}

We now can establish a good correspondence between $S$-maximal amenable subgroups of $\mf G$ and maximal amenable subgroups of $\mf G(K(S))$.

\begin{Cor}\label{cor:all MA in K(S) points are H(K(S))}
Every maximal amenable subgroup $\Lambda<\mf G(K(S))$ is of the form $\mf H(K(S))$ for some $S$-maximal amenable subgroup $\mf H<\mf G$.
\end{Cor}
\begin{proof}
 The Zariski closure of $\Lambda$ is $K$-defined \cite[Corollary AG.14.6]{BorelLinear2012}.
 By Proposition \ref{prop: poor mans tits alternative}, the Zariski closure of $\Lambda$ is $S$-amenable. By Proposition \ref{prop: every SA contained in SMA} such a group is contained in a $K$-defined $S$-maximal amenable subgroup $\mf H\leqslant\mf G$. Since $\mf H(K(S))\geqslant \Lambda$ is amenable by Lemma \ref{lem: amenable iff K(S) points amenable}, we must have that $\mf H(K(S))=\Lambda$. 
\end{proof}

The converse of the above corollary is not clear. The sitaution is much easier in the case of \textit{commensurably maximal amenable} subgroups.

Recall that two subgroups of a group are called commensurable if their intersection has finite index in both of them; commensurability is an equivalence relation on the set of subgroups of a group. Clearly all groups in the same commensurability class are either amenable or not. Thus, one can speak of amenable and non-amenable commensurability classes.

\begin{Def}\label{def: commensurably maximal amenable}
A maximal amenable subgroup $\Lambda\leqslant \Gamma$ is called \textbf{commensurably maximal amenable} if its commensurability class is maximal.
\end{Def}

Since increasing unions of amenable subgroups are amenable, every commensurably maximal amenable subgroup is contained in a maximal amenable subgroup as a subgroup of finite index.

\begin{Prop}\label{prop: every CMA is H(K(S))}
  A subgroup $\Lambda<\mf G(K(S))$ is commensurably maximal amenable if and only if there is an $S$-cocompactly maximal amenable subgroup $\mf H<\mf G$ such that $\Lambda=\mf H(K(S))$.
\end{Prop} 
\begin{proof}
Let us assume that $\Lambda$ is commensurably maximal amenable. 
By Corollary \ref{cor:all MA in K(S) points are H(K(S))}, $\Lambda=\mf H(K(S))$ for some $S$-maximal amenable group $\mf H$. Set $\mf U\coloneqq \mc R_u(\mf H)$. By Theorem \ref{thm:maximal-amenable-algebraic-subgroups} \ref{mainthm: S maximal amenable} $\mf H=\mf H^R\ltimes \mc R_u(\mf H)$ and $\mf H^R$ contains an $S$-cocompact torus $\mf T$. We can also assume that $\mf T$ has no $S$-cocompact subgroups; then $\mf T\ltimes \mf U$ has no $S$-cocompact subgroups either. If $\mf H$ is not $S$-cocompactly maximal amenable, by Proposition \ref{prop: every sma is contained scma} there is an $S$-cocompactly maximal amenable $\mf H'\leqslant\mf G$ such that $\mf H\cap \mf H'$ is $S$-cocompact in $\mf H$ and $\mf H'/(\mf H\cap \mf H')$ is a non-$S$-compact torus of $K$-rank $0$. Hence $\Lambda\cap \mf H'(K(S))$ has finite index in $\Lambda$ but infinite index in $\mf H'(K(S))$ which contradicts our hypothesis.

Conversely, if $\mf H$ is $S$-cocompactly maximal amenable, then by Theorem \ref{thm:maximal-amenable-algebraic-subgroups} \ref{mainthm: S cocompactly maximal amenable} $\Lambda$ is commensurable to a group of the form $\mf H_{\mf T}(K(S))$. Any of its finite index subgroups contains $\mf U$ and the $K$-anisotropic part $\mf T_a$ of $\mf T$ in its Zariski closure, see Remark \ref{rem:HT(K(S)) Zariski dense in HT}. Thus, any amenable overgroup $\Lambda'$ of a finite index subgroup of $\Lambda$ has to contain $\mf T_a\ltimes \mf U$ in its Zariski closure. In particular, if $\Lambda'$ is maximal amenable, by Corollary \ref{cor:all MA in K(S) points are H(K(S))} and Theorem \ref{thm:maximal-amenable-algebraic-subgroups} \ref{item2: normalizer is MA} we conclude that $\Lambda'$ is commensurable to $\Lambda$.
\end{proof}

\subsection{Singularity} 

The following definition is borrowed from \cite{BoutonnetMaximal2015}. 

\begin{Def}
    An amenable subgroup $\Lambda<\Gamma$ is \textbf{singular} if there is a continuous action $\Gamma\curvearrowright X$ on a compact space such that for any $\Lambda$-invariant probability measure $\mu$ and for all $\gamma\in \Gamma\setminus \Lambda$, the probability measures $\mu$ and $\gamma_*\mu$ are singular.
\end{Def}

Our aim is to prove that all commensurably maximal amenable groups of arithmetic groups are singular, Theorem \ref{ccma singular K}. Before proving it, let us remark that the converse is always true. 

\begin{Lemma}
 Assume that the amenable group $\Lambda$ is singular in $\Gamma$. Then $\Lambda$ is commensurably maximal amenable in $\Gamma$.
\end{Lemma}
\begin{proof}
    Consider the action $\Gamma\curvearrowright X$ as in the definition of singularity, let $\Lambda'\leqslant\Lambda$ be a finite index subgroup and take an amenable group $\Lambda''\geqslant \Lambda'$. Since $\Lambda''$ is amenable, there is a $\Lambda''$-invariant measure $\mu$ on $X$. Clearly $\mu$ is $\Lambda'$-invariant, hence its orbit $\Lambda\mu=\{\mu_0,\mu_1,\ldots,\mu_k\}$ is finite, where $\mu_0=\mu$. Set $\mu'\coloneqq 1/(k+1)\sum_i \mu_i$. Then $\mu'$ is $\Lambda$-invariant. Given any $\gamma\in \Lambda''$, remark that the measures $\mu'$ and $\gamma_*\mu'$ are not singular: they are both the convex sum of $\mu$ with another probability measure. Therefore $\Lambda''\leqslant\Lambda$ and hence $\Lambda$ is commensurably maximal amenable.  
\end{proof}

Our main theorem is the converse of the above lemma for arithmetic lattices.

\begin{Thm}\label{ccma singular K}
  Let $K$ be a global field of characteristic $0$ and let $S\subseteq \ms R$ be a set of places of $K$ satisfying $\ms R_\infty\setminus \ms T(\mf G)\subseteq S\subseteq \ms R$. Let $\mf G$ be a connected $K$-group and $\mf H$ be an $S$-cocompactly maximal amenable subgroup. Then $\mf H(K(S))$ is a singular subgroup of $\mf G(K(S))$.
\end{Thm}  
\begin{proof}
As described in Theorem \ref{thm:maximal-amenable-algebraic-subgroups} \ref{mainthm: S cocompactly maximal amenable}, $\mf H$ is of the form $\mf H = \mf H^R \ltimes \mf U$, where $\mf H^R = \mc N_{\mf L}(\mf T)$ is the normaliser of an $S$-ample torus $\mf T$ in $\mf L$, where $\mf L$ is the centralizer of the maximal $K$-split subtorus of $\mf T$. We denote by $\mf T_a$ the maximal $K$-anisotropic subtorus of $\mf T$; we also retain the notation $\mf P = \mc N_{\mf G}(\mf U)$. For each $v\in S$, we denote by $\mf T_{v}$ the maximal $K_v$-split subtorus of $\mf T$ and by $\mf T_{a,v}$ the maximal $K_v$-split subtorus of $\mf T_a$. 

Choose a minimal $K_v$-parabolic subgroup $\mf Q_v\leqslant \mf P$ containing $\mf T_v$ and set
\[
X_v \coloneqq \mf G(K_v)/\mf Q_v(K_v),\quad X\coloneqq\prod_{v\in S} X_v
\]
eqipped with the diagonal action of $\mf G(K)$.

Let $\mu$ be an $\mf H(K(S))$-invariant probability measure on $X$. Take an arbitrary $\gamma\in \mf G(K)$ and assume that $\gamma_*\mu$ and $\mu$ are not singular; we will show that $\gamma\in\mf H(K)$. Let $\mu_v$ be the projection of $\mu$ on $X_v$; clearly, $\mu_v$ is $\mf H(K(S))$-invariant. Then for every $v\in S$ the measures $\gamma_*\mu_v$ and $\mu_v$ are not singular by the following claim.

\begin{Claim}
    Let $\nu$ and $\nu'$ be probability measures on $X$ and denote as above by $\nu_v$ and $\nu_v'$ their projections on $X_v$. Then if $\nu$ and $\nu'$ are not singular, then for all $v\in S$ the measures $\nu_v$ and $\nu'_v$ are not singular.
\end{Claim}
\begin{proof}[Proof of Claim]
    If there is a $v\in S$ such that $\nu_v$ is singular with respect to $\nu'_{v}$, then there is a commonly measurable subset $A_v\subseteq X_v$ such that $\nu_{v}(A_v) = 0$ and $\nu'_v(A_v) = 1$. But then 
    \[
        \nu\left(A_v\times \prod_{w\neq v} X_w\right) = 0,\quad \nu'\left(A_v\times \prod_{w\neq v} X_w\right) = 1
    \]
    which proves that $\nu$ and $\nu'$ are singular.
\end{proof}
Fix an arbitrary $v\in S$.
\begin{Claim}
    The stabilizer $G_{\mu_v}\leqslant \mf G(K_v)$ of the measure $\mu_v$ contains $\mf T_{a,v}(K_v)\ltimes \mf U(K_v)$. 
\end{Claim}
\begin{proof}[Proof of Claim]
By \cite[Theorem 3.2.4]{ZimmerErgodic1984} (see also \cite[Corollary 1.8]{BaderAlmost2017}), $G_{\mu_v}$ contains a normal subgroup $G_{\mu_v}'$ of finite index which is almost $K_v$-algebraic. Let us denote by $\mf N$ a normal $K_v$-defined $K_v$-cocompact subgroup of $G'_{\mu_v}$. By Remark \ref{rem:HT(K(S)) Zariski dense in HT}, the group $(\mf T_a\ltimes \mf U)(K(S))$ is Zariski dense in $\mf T_a\ltimes \mf U$. Since a finite index subgroup of $(\mf T_{a}\ltimes \mf U)(K(S))$ is contained in $G'_{\mu_v}$ and therefore normalizes $\mf N(K_v)$, its Zariski closure $\mf T_a\ltimes \mf U$ normalizes $\mf N$, that is, $\mf T_a\ltimes\mf U \leqslant \mc N_{\mf G}(\mf N)$. Set 
\[\mf B\coloneqq (\mf T_a\ltimes \mf U)/(\mf N\cap (\mf T_a\ltimes \mf U)).\]

By assumption on $\mf N$, the quotient $G'_{\mu_v}/\mf N(K_v)$ is compact.
Therefore the image $\Lambda$ of $(\mf T_a\ltimes \mf U)(K(S))$ in $\mf G(K_v)/\mf N(K_v)$ is a pre-compact subgroup of $\mf B(K_v)$.

By \cite[Theorem I.3.2.4, Remark I.3.2.6]{MargulisDiscrete1991}, there is a compact subset $C\leqslant (\mf T_a\ltimes\mf U)(K_v)$ such that 
\[(\mf T_a\ltimes \mf U)(K(S))\cdot C=(\mf T_a\ltimes \mf U)(K_v).\]
Clearly $\Lambda [C]_{\mf N}=\mf B(K_v)$ and 
since the image $[C]_{\mf N}$ of the compact set $C$ is compact, and $\Lambda$ is precompact in $\mf B(K_v)$, we obtain that $\mf B(K_v)$ is compact. 

This implies that $\mf N\cap (\mf T_a\ltimes\mf U)$ is a $K_v$-defined $K_v$-cocompact subgroup of $\mf T_a\ltimes\mf U$. Such a subgroup has to contain $\mf T_{a,v}\ltimes \mf U$, since the latter has no $K_v$-defined $K_v$-cocompact subgroups.  
\end{proof} 

By \cite[Theorem 3.1.3, Proposition 2.1.12]{ZimmerErgodic1984}, the action of $(\mf T_{a,v}\ltimes \mf U)(K_v)$ is algebraic and hence smooth on $X_v$. Therefore any $(\mf T_{a,v}\ltimes\mf U)(K_v)$-invariant probability measure has to be contained in $(\mf T_{a,v}\ltimes\mf U)(K_v)$-compact orbits. Such a group does not have any proper $K_v$-defined $K_v$-cocompact subgroups, so that $\mu_v$ is supported on the set $F$ of fixed points of $(\mf T_{a,v}\ltimes \mf U)(K_v)$.  

\begin{Claim}\label{claim: g in K_v}
  If $[g]\in F$, then $g\in \mf P(K_v)$. In particular, every $[g]\in F$ is fixed by $(\mf T_v\ltimes\mf U)(K_v)$.
\end{Claim}
\begin{proof}[Proof of Claim]
Note that $g\in\mf G(K_v)$ is such that $[g]\in F$ if and only if $(\mf T_{a,v}\ltimes \mf U)(K_v)\subseteq g\mf Q_v(K_v)g^{-1}$. Now, observing that $\mf Q_v \leqslant \mf P$, we get the following: if $\mf U \leqslant g\mf Q_v g^{-1} \leqslant g\mf P g^{-1}$, then $\mf U \leqslant g\mf P g^{-1} \cap \mf P$. Now, \cite[Corollaire 4.5]{BorelGroupes1965} ensures that $g\mf P g^{-1} = \mf P$ and as $\mf P$ is self-normalizing \cite[Proposition 4.4]{BorelGroupes1965}, $g\in \mf P$. 
Therefore if $g\in\mf P(K_v)$, we can write $g=\ell u$ where $\ell\in\mf L(K_v)$ and $u\in \mf U(K_v)\leqslant\mc R_u(\mf Q_v)(K_v)$, hence $[g]=[\ell]$. Since the $K$-split subtorus of $\mf T$ commutes with $\ell$ and it is contained in $\mf Q_v$, it follows that every point in $F$ is fixed by $\mf T_v\ltimes\mf U$.
\end{proof}

It follows from the claim that the measure $\mu_v$ is supported on the set $F$ of fixed points of $(\mf T_v\ltimes \mf U)(K_v)$. As $\gamma_*\mu_v$ and $\mu_v$ are not singular, we must have that $\gamma F\cap F$ is not empty. In particular, there is a $g\in \mf G(K_v)$ such that $[g]\in F$ and $[\gamma g]\in F$. By Claim \ref{claim: g in K_v}, we have $g,\gamma g\in \mf P(K_v)$, so $\gamma\in \mf P(K)$. Decomposing $\gamma = \gamma_{\ell}\gamma_u$ with $\gamma_\ell\in \mf L(K)$ and $\gamma_u\in\mf U(K)$ and noticing that $\mf U(K)$ fixes $F$ pointwise, we see that $\gamma_*\mu_v =(\gamma_\ell)_*\mu_v\eqqcolon \mu'_v$. We let $\mf T'\coloneqq  \gamma_\ell \mf T \gamma_\ell^{-1}\leqslant \mf L$ and use the subscripts for the conjugates of corresponding subtori of $\mf T$.

Let us denote by $H_{\mu'_v}\leqslant \mf T'(K_v)\ltimes \mf U(K_v)$ the stabilizer of the measure $\mu'_v$.

\begin{Claim}
    The group $H_{\mu'_v}$ is Zariski dense and $K_v$-cocompact in $\mf T'\ltimes \mf U$. 
\end{Claim}
\begin{proof}[Proof of Claim]
  The support of $\mu_v'$ is pointwise fixed by $\mf T_v' \ltimes \mf U$ which is $K_v$-cocompact in $\mf T'\ltimes\mf U$. Since $\mf T'$ is $S$-ample, $(\mf T'_a\ltimes \mf U)(K(S))$ is Zariski dense in $\mf T'_a\ltimes \mf U$, see Remark \ref{rem:HT(K(S)) Zariski dense in HT}. Since $\mf T'_v$ and $\mf T'_a$ generate $\mf T'$, the claim follows.
\end{proof}

Note that the action of $H_{\mu'_v}$ on the support of $\mu'_v$ is smooth: its cocompact normal subgroup $\mf T'_v\ltimes \mf U(K_v)$ acts trivially and therefore it reduces to an action of a compact group \cite[Corollary 2.1.13]{ZimmerErgodic1984}. Take the $H_{\mu'_v}$-ergodic decomposition \[\mu_v'=\int_Y \nu_y d\eta.\]
Each $\nu_y$ is supported on a single $H_{\mu'_v}$-orbit,
say $H_{\mu'_v}[g_y]$ for some $g_y\in \mf G(K_v)$; 
let us denote by $\mf V_y$ the Zariski closure of this orbit. Since $H_{\mu'_v}\leqslant \mf T'(K_v)\ltimes \mf U(K_v)$ is Zariski dense by the above claim, we also have that $(\mf T'\ltimes \mf U)[g_y]\subseteq \mf V_y$, so $\mf V_y$ is also the Zariski closure of $(\mf T'\ltimes \mf U)[g_y]$. Note that since $\mf T'\ltimes \mf U$ acts on $\mf V_y$, it permutes the irreducible components and since $\mf T'\ltimes \mf U$ is Zariski-connected, it follows that $\mf V_y$ is irreducible.  

 We denote by $\mf M_v$ the centralizer of $\mf T_v$. Since $\mf L$ is the centralizer of the maximal $K$-split torus of $\mf T$ and the latter is contained in $\mf T_v$, we have that $\mf M_v\leqslant \mf L\leqslant \mf P$. Hence $\mf M_v$ normalizes $\mf U$, which altogether implies that $\mf M_v$ normalizes $\mf T_v\ltimes \mf U$, therefore $\mf M_v(K_v)$ preserves its fixed point set $F$.

\begin{Claim}\label{claim: measure V_y is 0}
For every $[g]\in F$,
  the orbit $\mf M_v[g]$ is an irreducible projective variety. If $\mf V_y$ is not contained in $\mf M_v[g]$, then $\nu_y(\mf M_v(K_v)[g])=0$.
\end{Claim}
\begin{proof}[Proof of Claim]
 To show that the orbit of $[g]$ is an irreducible projective subvariety, observe that $\mf T_v\leqslant g\mf Q_v g^{-1}$ and let $\mf R_v$ be a parabolic subgroup containing $g\mf Q_v g^{-1}$ whose Levi subgroup is $\mf M_v$. By \cite[Proposition 4.4 c)]{BorelGroupes1965}, every element $r\in \mf R_v(K_v)$ has a decomposition $r = mu$, where $m\in \mf M_v(K_v)$ and $u\in \mc R_u(\mf Q_v)(K_v)$ which clearly implies that $\mf R_v(K_v)[g]=\mf M_v(K_v)[g]$. In particular, since $\mf R_v/g\mf Q_v g^{-1}$ is an irreducible projective $K_v$-variety, so is $\mf M_v(K_v)[g]$.

 Let now $\mf W\subseteq \mf V_y$ be an irreducible subvariety of minimal dimension with the property that $\nu_y(\mf W)$ is positive. Since $\nu_y$ is $H_{\mu'_v}$-invariant, for every $h\in H_{\mu'_v}$ we have that $\nu_y(\mf W)=\nu_y(h\mf W)$. Also observe that $g\mf W$ is an irreducible variety and therefore $\mf W\cap h\mf W$ has smaller dimension than $\mf W$ and hence is a $\nu_y$-null set unless $\mf W=h\mf W$. By ergodicity, we get finitely many translates $h_1\mf W,\ldots,h_n\mf W$ whose union is of measure $1$. Then $\mf V_y = \bigcup_{i=1}^n h_i\mf W$, and by irreducibility it follows that $\mf V_y = \mf W$. If $\mf V_y$ is not contained in $\mf M_v[g]$, their intersection has strictly lower dimension, hence is of measure zero. 
\end{proof}

\begin{Claim}\label{claim: finitely many Mv-orbits}
  There are only finitely many $\mf M_v(K_v)$ distinct orbits in $F$.
\end{Claim}
\begin{proof}[Proof of Claim]
  Fix an arbitrary point $[g]\in F$. Let $\mf S\leqslant \mf Q_v$ be a maximal $K_v$-split torus containing $\mf T_v$. As we have already remarked, $g^{-1}\mf T_v(K_v) g\subseteq \mf Q_v(K_v)$. Since $\mf T_v$ is $K_v$-split, there is a unipotent element $u\in \mc R_u(\mf Q_v)(K_v)$ such that $u^{-1}g^{-1}\mf T_v(K_v)gu\subseteq \mf S(K_v)$. Also note that $[g]=[gu]$. Therefore every point in $F$ is represented by an element of the \textit{transporter} of $\mf T_v$ into $\mf S$. By \cite[Proof of Proposition 13.3.1]{SpringerLinear1998}, this transporter is the union of finitely many $\mf M_v$-cosets.
\end{proof}

Since $\mu_v$ and $\mu'_v$ are not singular, we have that $\mu'_v$-measure of the support of $\mu_v$ cannot be zero. By Claim \ref{claim: finitely many Mv-orbits}, the support of $\mu_v$ is a finite union of $\mf M_v(K_v)$-orbits, so there is at least one $\mf M_v(K_v)$-orbit which has positive measure with respect to both $\mu_v$ and $\mu'_v$. Since $\mu_v$ is supported on $F$, it has to be an $\mf M_v(K_v)$-orbit of an element $[g]\in F$.
Claim \ref{claim: measure V_y is 0} implies that there is $y$ such that $(\mf T'\ltimes \mf U)(K_v)[g_y]\subseteq \mf M_v(K_v)[g]\subseteq F$. As $[g_y]\in \mf M_v(K_v)[g]$, we also have that $[g_y]\in F$ and $(\mf T'\ltimes \mf U)[g_y]\subseteq \mf M_v[g_y]$.
As in the proof of Claim \ref{claim: measure V_y is 0}, we consider the parabolic subgroup $\mf R_v$ containing $g_y\mf Q_v g_y^{-1}$ whose Levi subgroup is $\mf M_v$; as already observed, $\mf R_v(K_v)[g_y]=\mf M_v(K_v)[g_y]$. 

\begin{Claim}
    We have $\mf T'\cap \mf T\supseteq \mf T_v$; in particular, $\mf T_v = \mf T'_v$. 
\end{Claim}
\begin{proof}[Proof of Claim]
    As observed before, $\mf T'[g_y]\subseteq \mf M_v[g_y]$. Now since $\mf R_v\geqslant g_y\mf Q_vg_y^{-1}$ and $\mf M_v[g_y]=\mf R_v[g_y]$, we have $\mf T'\leqslant \mf R_v$. Let $\mf N\leqslant \mf R_v$ be the $K$-group generated by $\mf T'$ and $\mf T$. Since $\mf T$ is $S$-ample in $\mf L$, $\mf T$ cannot normalize any unipotent $K$-subgroup in $\mf L$. Therefore $\mf N$ is reductive. Then there is an element $u\in \mc R_u(\mf R_v)(K_v)$ such that $\mf N\leqslant u\mf M_v u^{-1}$. Since $\mf T$ and $\mf T'$ are $S$-ample, they have to contain the the $K_v$-split part of the central torus of $u\mf M_v u^{-1}$. By rank considerations the latter then has to be equal to both $\mf T_v$ and $\mf T'_v$. 
\end{proof}

We now deduce that $\gamma_\ell$ normalizes $\mf T_v$ for every $v\in S$. As $\mf T$ is $S$-ample, it coincides with the minimal $K$-defined torus containing all $\mf T_{v}$, and thus we see that $\gamma_\ell\in \mc N_{\mf L}(\mf T)(K)$. Therefore $\gamma\in (\mc N_{\mf L}(\mf T)\ltimes \mf U)(K) = \mf H(K)$. This finishes the proof.
\end{proof}

\section{Some examples}
Constructing explicit examples of commensurably maximal amenable subgroups of arithmetic groups in classic algebraic groups can be done using number-theoretic input and computer algebra assistance. Implementing algorithms which find $S$-ample tori and giving explicit generators of commensurably maximal amenable subgroups in general arithmetic groups will be subject of forthcoming independent work. In this section we will restrict our attention to giving some easy examples of commensurably maximal amenable subgroups in $\mf{SL}_n(\mb Z)$ and $\mf{SL}_n(\mb Z[1/5])$.

Constructing maximal tori in $\mf{GL}_n(\mb Q)$ amount to the study of étale algebras. Recall that an étale algebra $E$ over $\mb Q$ is a finite product of algebraic extensions of $\mb Q$:
\[
E = \prod_{k=1}^\ell E_k,
\]
where $E_k$ is an algebraic extension of $\mb Q$ of degree $n_k\coloneqq  [E_k:\mb Q]$. The degree of $E$ is $n\coloneqq [E:\mb Q] = \dim_{\mb Q} E = \sum_{k=1}^\ell n_k$, and fixing a $\mb Q$-basis of $E$ defines a matrix representation
\[
\pi\colon E\to \mb M_{d}(\mb Q)
\]
coming from the regular representation of $E$.

Then there is a unique maximal torus $\mf T < \mf {GL}_n$ such that
\[
\mf T(\mb Q) = \pi(E^\times)
\]
and every $\mb Q$-defined maximal torus in $\mf{GL}_n$ is of this form. Indeed, the algebra generated by $\mf T(\mb Q)$ inside $\mf M_n(\mb Q)$ is easily seen to be étale since $\mf T$ is $\overline{\mb Q}$-diagonalizable.

We now want to investigate the $\mb Z$-points of such a torus. In general, these depend on the choice of the matrix representation $\pi$ or, equivalently, on the choice of a $\mb Q$-basis of $E$.

Given a matrix representation $\pi$ as above, we obtain an order of $E$
\[
\mc O\coloneqq \{ e\in E \mid \pi(e)\in \mb M_n(\mb Z)\}
\]
which clearly satisfies
\[
\mf T(\mb Z) = \pi(\mc O^\times).
\]

Conversely, choosing an order $\mc O\subset E$ and considering the matrix representation $\pi$ associated to a $\mb Z$-basis of $\mc O$, we have
\[
\{e\in E\mid \pi(e) \in \mf T(\mb Z)\} = \mc O^\times
\]
Choosing a different $\mb Z$-basis of $\mc O$ yields a torus $\mf T'$ such that $\mf T'(\mb Z)$ and $\mf T(\mb Z)$ are conjugate inside $\mf{GL}_n(\mb Z)$. Therefore, fixing a maximal $\mb Q$-torus inside $\mf {GL}_n$ up to conjugacy by $\mf{GL}_n(\mb Z)$ amounts to fixing an étale algebra $E$ together with an order $\mc O$. The $\mb Q$-rank of such a torus $\mf T$ is equal to $\rk_{\mb Q}\mf T=\ell$.

Every automorphism $\sigma\in \Aut(\mc O)\leqslant \Aut(E)$ yields a matrix $\pi(\sigma)\in \mf{GL}_n(\mb Z)$ which normalizes $\pi(E)$ and therefore lies in $\mc N_{\mf GL_n}(\mf T)(\mb Z)$. Conversely, if $g\in \mc N_{\mf GL_n}(\mf T)(\mb Z)$, then conjugation by $g$ defines an automorphism of $\pi(\mc O)\cong \mc O$, and therefore $g = \pi(\sigma)\cdot \pi(t)$ for some $t\in \mf T(\mb Z)$. 

Let us begin with the case when $S$ consists only of the real valuation. A maximal $\mb Q$-defined and $\mb Q$-irreducible torus $\mf T$ of $\mf{SL}_n$ is $S$-ample whenever it is not $\mb R$-compact. Note that the latter can happen only for $n=2$. Let us give examples of such tori and the corresponding maximal amenable subgroups of $\mf{SL}_n(\mb Z)$. In our examples the order $\mc O$ will always be the maximal order in $E$, so that the generators of $\mf T(\mb Z)$ are the generators of the unit group of $E$; they were computed using the \textsc{Sage} computer algebra system.

\begin{Ex}\label{ex: torus is RxS1}
Taking the torus associated to the number field $E = \mb Q[x]/\langle x^3+x-1\rangle$, we obtain a commensurably maximal amenable subgroup of $\mf{SL}_3(\mb Z)$ isomorphic to $\mb Z$ generated by
\[
g = \begin{pmatrix}
0 & 0 & 1\\
1 & 0 & -1\\
0 & 1 & 0
\end{pmatrix}.
\]
There are no integral elements in the normalizer because the automorphism group of $E$ is trivial.
\end{Ex}

\begin{Ex}
Consider the torus in $\mf{SL}_4$ associated to the number field $E = \mb Q[x]/\langle x^4 - 8x^3 + 20x^2-16x + 1\rangle$. Its integer points are generated by the central element $-\mf 1$ and
\[
g_1=
\begin{pmatrix}
0  & -1  & -2 & -4 \\   
8  & 16  & 31 & 62 \\   
-6 & -12 & -24& -49 \\   
1  & 2   & 4  & 8 \\
\end{pmatrix},
g_2=
\begin{pmatrix}
-1 & -1 & -2 & -5\\
 9 & 15 & 31 & 78\\
-6 &-11 &-25 &-69\\
 1 &  2 &  5 & 15\\
\end{pmatrix},
g_3=
\begin{pmatrix}
-4 & -1 & -1 & -1\\
13 & 12 & 15 & 15\\
-7 & -7 & -8 & -5\\
 1 & 1  & 1  & 0
\end{pmatrix}.
\]
Here the Galois group $\Gal(K/\mb Q) \cong \mb Z/2\times\mb Z/2$ generated by
\[
w = 
\begin{pmatrix}
  1 &  0 & 0 &  0\\
-16 & -1 &-4 &-16\\
 16 &  0 & 1 & 8\\
 -4 &  0 & 0 &-1
\end{pmatrix},
w' =
\begin{pmatrix}
 1  & 0 &  0 &  0 \\
-14 &  4 &  10 & 25\\
 12 & -4 & -9 & -20\\
-2  & 1  & 2  & 4
\end{pmatrix}.
\]

The group generated by $-\mf 1,g_1,g_2,g_3,w$ and $w'$ is commensurably maximal amenable.
\end{Ex}

By considering the torus of Example \ref{ex: torus is RxS1} inside $\mf{SL}_4$, we obtain an example of a maximal amenable group with unipotent elements. 

\begin{Ex}
  We let $\tilde g\in \mf{SL}_4(\mb Z)$ be the block diagonal matrix 
  \[
  \tilde g=
  \begin{pmatrix}
  g&0\\
  0&1
  \end{pmatrix}
  \]
  where $g$ is as in Example \ref{ex: torus is RxS1}. If we denote by $E_{i,j}\in\mf{SL}_4(\mb Z)$ the standard elementary matrix, then one has that the group generated by \[\tilde g, -\mf 1, E_{1,4}, E_{2,4}, E_{3,4}\] is commensurably maximal amenable. 
\end{Ex}

Let us finally give an example where $S$ consists of $2$ valuations.

\begin{Ex}
Consider the maximal $\mb Q$-torus $\mf T$ in $\mf{SL}_2$ associated with the number field $E = \mb Q[i]$. Then we have $\rk_{\mb R}\mf T = 0$ (so it is not ample for $S=\{\infty\}$) and $\rk_{\mb Q_5} = 1$. In particular, $\mf T$ is $\{\infty,5\}$-ample and $\mf T(\mb Z[1/5])$ is generated by 
\[
i =
\begin{pmatrix}
0 & -1  \\ 1 & 0 
\end{pmatrix},\quad g =\begin{pmatrix}
4/5 & -3/5\\
 3/5 & 4/5\\
\end{pmatrix}
\]
as can be verified by computing $S$-units of norm $1$ in $E$. This group is commensurably maximal amenable.
\end{Ex}

\end{document}